\documentclass[11pt]{article}
\usepackage{amsfonts}
\usepackage{dsfont}
\usepackage{amsfonts}
\usepackage{bbm}
\usepackage{mathrsfs}
\usepackage{color}
\usepackage{hyperref}
\usepackage{mathrsfs}
\usepackage{amsmath}
\usepackage{enumitem}
\allowdisplaybreaks[4]
\usepackage{amsfonts}
\usepackage{amsthm}
\usepackage{amssymb, amsmath, cite}
\usepackage{bm}
\usepackage{amsthm}
\usepackage{verbatim}
\usepackage{graphicx}
\usepackage{color}
\usepackage{amsfonts}
\usepackage{dsfont}
\usepackage{amsfonts}
\usepackage{bbm}
\usepackage{mathrsfs}
\usepackage{color}
\usepackage{graphicx,amsmath,amsfonts,amssymb,amscd,amsthm}
\usepackage{hyperref}
\usepackage[numbers,sort&compress]{natbib}
\usepackage{enumitem}
\setenumerate[1]{itemsep=1pt,partopsep=0pt,parsep=\parskip,topsep=1pt}
\setitemize[1]{itemsep=1pt,partopsep=0pt,parsep=\parskip,topsep=1pt}

\newtheorem{Theorem}{Theorem}[section]
\newtheorem{Lemma}{Lemma}[section]

\newtheorem{Remark}{Remark}[section]
\newtheorem{Definition}{Definition}[section]

\usepackage{graphicx,amsmath,amsfonts,amssymb,amscd,amsthm}
\usepackage[numbers,sort&compress]{natbib}
\usepackage{array}
\usepackage{booktabs}
\usepackage{multirow}
\usepackage{makecell}
\usepackage{abstract}
\hypersetup{colorlinks,linkcolor=blue,citecolor=blue}
\newtheoremstyle{kai}
{3pt} {3pt} {} {} {\bfseries} {.} {.5em} {}
\def\EquationsBySection{\def\theequation
{\thesection.\arabic{equation}}
\@addtoreset{equation}{section}}

\newcommand\old[1]{}

\newcommand{\Bp}{\begin{proof}}
 \newcommand{\Ep}{\end{proof}}
\renewcommand{\theequation}{\thesection.\arabic{equation}}

\numberwithin{equation}{section}
\addcontentsline{toc}{section}{Acknowledgement}

\newcommand{\newcom}{\newcommand}
\newcom{\R}{\mathbb R}
\newcom{\N}{\mathbb N}
\newcom{\e}{\varepsilon}

\newcom{\al}{\alpha}
\newcom{\be}{\beta}
\newcom{\del}{\delta}

\newcom{\ga}{\gamma}
\newcom{\Ga}{\Gamma}

\newcom{\Lam}{\Lambda}
\newcom{\lam}{\lambda}

\newcom{\Om}{\Omega}
\newcom{\om}{\omega}

\newcom{\Si}{\Sigma}
\newcom{\si}{\sigma}
\newcom{\s}{\varsigma}

\newcom{\tht}{\theta}
\newcom{\dtri}{\nabla}
\newcom{\tri}{\triangle}
\newcom{\mf}{\,\,\,\,\mathrm{for\,\,all} \,\,}
\newcom{\qmf}{\qquad\mathrm{for\,\,all}~~ }

\newcom{\no}{\nonumber}
\newcom{\f}{\frac}
\newcom{\na}{\nabla}
\newcom{\Del}{\Delta}
\newcom{\ep}{\epsilon}

\newcom{\p}{\partial}

\newcom{\uep}{u_{\e}}
\newcom{\vep}{v_{\e}}
\newcom{\nep}{n_{\epsilon}}
\newcom{\cep}{{c}_{\epsilon}}
\newcom{\whc}{\widehat{c}}
\newcom{\ocep}{\widehat{c}_{\epsilon}}
\newcom{\cp}{{c}_{\epsilon}}
\newcom{\iom}{\int_{\Omega}}
\newcom{\wcep}{\widetilde{c}_{\e}}
\newcom{\ce}{c^*}

\newcom{\beq}{\begin{equation}}
\newcom{\eeq}{\end{equation}}
\newcom{\beno}{\begin{eqnarray*}}
\newcom{\eeno}{\end{eqnarray*}}

\newcommand{\bea}{\begin{eqnarray}\label}
\newcommand{\eea}{\end{eqnarray}}
\newcommand{\bas}{\begin{eqnarray*}}
\newcommand{\eas}{\end{eqnarray*}}

\begin{document}
\title{\bf The asymptotic behavior of solutions to a doubly degenerate chemotaxis-consumption system in the two-dimensional setting}

\author{{\sc Duan Wu\thanks{E-mail: duan\underline{~}wu@126.com}}\\
{\small Institut f\"ur Mathematik, Universit\"at Paderborn, 33098 Paderborn, Germany}}
\date{}
\maketitle
\begin{abstract}
 The present work proceeds to consider the convergence of the solutions to the following doubly degenerate chemotaxis-consumption system
\begin{align*}
\left\{
\begin{array}{r@{\,}l@{\quad}l@{\,}c}
&u_{t}=\nabla\cdot\big(u^{m-1}v\nabla v\big)-\nabla\cdot\big(f(u)v\nabla v\big)+\ell uv,\\
&v_{t}=\Delta v-uv,
\end{array}\right.
\end{align*}
under no-flux boundary conditions in a smoothly bounded convex domain $\Omega\subset \R^2$, where the nonnegative function $f\in C^1([0,\infty))$ is asked to satisfy $f(s)\le C_fs^{\al}$ with $\al, C_f>0$ for all $s\ge 1$.

The global existence of weak solutions or classical solutions to the above system has been established in both one- and two-dimensional bounded convex domains in previous works. However, the results concerning the large time behavior are still constrained to one dimension due to the lack of a Harnack-type inequality in the two-dimensional case. In this note, we complement this result by using the Moser iteration technique and building a new Harnack-type inequality.\\
   \\
\textbf{Keywords}: degenerate diffusion; chemotaxis; asymptotic stability\\
\textbf{AMS (2020) Subject Classification}: {35K65, 92C17, 35B40}
\end{abstract}

\section{Introduction}
Investigating the long-term behavior of bounded solutions {to} a parabolic system can serve as one of the criteria for assessing whether or not the system is capable of precisely describing the emergence and evolution of structures in some specific biological {situations}. {For instance, the concordance between the nontrivial stabilization characteristics shown in \cite{winkler2022cvpde} and the experimental observations reported in \cite{fujikawa1992physica,m-fujikawa1992physica,ohgiwari-m-m1992physica} rigorously indicates that, nonlinear reaction-diffusion systems involving cross-degeneracy proposed in \cite{kawasaki-mmus1997jtb} are more suitable for modeling the bacterial patterning of \textit{Bacillus subtilis} grown on thin agar plates than the general linear non-degenerate systems.}

The main purpose of our work is to consider the asymptotic behavior of the chemotaxis-consumption system
\begin{align}\label{s0}
\left\{
\begin{array}{r@{\,}l@{\quad}l@{\,}c}
&u_{t}=\nabla\cdot\left(D(u,v)\nabla u\right)-\nabla\cdot\left(S(u,v)\nabla v\right)+f(u,v),\\
&v_{t}=\Delta v-uv,
\end{array}\right.
\end{align}
that arises in mathematical biology as a model for describing the evolutionary processes of the species {\textit{Bacillus subtilis}} living in poor nutrient environments {(\cite{Leyva.etal2013physica})}. At the level of mathematical analysis for this model, there is some literature concentrating on the fundamental global solvability for different choices of $S(u,v)$, $D(u,v)$ and $f(u,v)$. When considering the nutrient-induced proliferation case (i.e., $f(u,v)=uv$), Winkler demonstrated that the one-dimensional version of \eqref{s0} admits global weak solutions in the context of $S(u,v)=uv$, $D(u,v)=u^2v$ in \cite{winkler2021tams}, which was subsequently followed by \cite{li-winkler2022cpaa} removing an integrability condition to obtain similar results; alternatively, global solvability in a two-dimensional setting can also be achieved by assuming smallness assumption on initial data, as shown in \cite{winkler2022na}. Also under the circumstances $D(u,v)=f(u,v)=uv$, the existence of global weak solutions holds {in} the taxis-free case (i.e., $S(u,v)=0$) in any dimension (\cite{winkler2022cvpde}); if chemotactic migration is taken into account of the form $S(u,v)=u^{\al}v$, Li in \cite{li2022jde} proved that \eqref{s0} possesses global weak solutions in two-dimensional domains with $1<\al<\f{3}{2}$ and in three-dimensional domains with $\f{7}{6}<\al<\f{13}{9}$, which {was} extended by \cite{winkler2024jde} and \cite{wu2024preprint}. In the presence of logistic source terms $f(u,v)=\rho u-\mu u^{\kappa}$, the global solvability was shown in \cite{pan2024norwa} for the case when $\kappa>\f{n+2}{2}$.

At the stage of research on asymptotic features, the available results from mathematically rigorous studies appear to be restricted to the one-dimensional version of \eqref{s0} with $D(u,v)=f(u,v)=uv$ and $S(u,v)=u^2v$ (\cite{winkler2021tams}), the striking peculiarity of which {consists in the presence of an uncountable family of spatially heterogeneous steady states significantly unlike the common case}. In the two-dimensional setting, only numerical simulations have been performed in \cite{Leyva.etal2013physica}, which coincides with the observations gained experimentally in \cite{ben-jacob.etal2000} and \cite{golding.etal1998}. To the best of our knowledge, however, {describing the role of steady states in the dynamics of \eqref{s0} by means of mathematical analysis in higher dimensions is} still lacking. The present work aims to address this gap.

More precisely, we continue our investigations related to the solutions of the initial-boundary value system
\begin{align}\label{s}
\left\{
\begin{array}{r@{\,}l@{\quad}l@{\,}c}
&u_{t}=\nabla\cdot\left(u^{m-1}v\nabla u\right)-\nabla\cdot\left(f(u)v\nabla v\right)+\ell uv,&x\in\Omega,&t>0,\\
&v_{t}=\Delta v-uv, &x\in\Omega,&t>0,\\
&\left(u^{m-1}v\nabla u-f(u)v\nabla v\right)\cdot\nu=\nabla v\cdot\nu=0,& x\in\partial\Omega,&t>0,\\
&u(x,0)=u_0(x),~~v(x,0)=v_0(x),&x\in\Omega,
\end{array}\right.
\end{align}
in a smoothly bounded convex domain $\Om\subset \R^2$, with $1 \le m<4$, $\ell\ge 0$, where the nonnegative function $f\in C^1([0,\infty))$ {is} assumed to satisfy
\beq\label{f1}
f(u)\le C_fu(u+1)^{\al-1} \qmf u\ge0
\eeq
or
\beq\label{f2}
f(u)\le C_fu^{\al} \qmf u\ge0
\eeq
with $C_f>0$ and $\al>0$. In addition, the initial data {are} throughout supposed to satisfy
\begin{align}\label{indata}
\left\{
\begin{array}{r@{\,}l@{\quad}l@{\,}c}
&u_0\in W^{1,\infty}(\Om){\mathrm{~is~nonnegative~with~}} u_0\not\equiv 0\qquad{\mathrm{and}}\\
&v_0\in W^{1,\infty}(\Om) {\mathrm{~is~positive~in~}}\overline{\Om}.
\end{array}\right.
\end{align}
{In} such frameworks, the global existence of weak solutions for $1\le m<3$ and classical solutions for $3\le m<4$ has been studied in \cite{wu2024preprint}. This paper will firstly illustrate that the solutions obtained previously actually exhibit higher regularities by employing a Moser iteration argument. These results are beneficial for establishing {a certain} Harnack-type inequality in the next step, which is crucial for {describing} the asymptotic behavior.

The following statement is with respect to the definition of weak solutions used in this paper.
\begin{Definition}\label{de}
Let $\Omega\subset \mathbb{R}^2$ be a bounded domain with smooth boundary, $m\ge 1$ and $\ell\ge 0$. Suppose that $f$ satisfies \eqref{f1} or \eqref{f2} with $C_f>0$ and $\al>0$, and that $u_0\in L^1(\Om)$ and $v_0\in L^1(\Om)$ are nonnegative. Then we call that a pair $(u,v)$ of nonnegative functions
\beno
\left\{
\begin{split}
&u\in L_{loc}^{1}(\overline\Om\times[0,\infty))\qquad{\rm{and}}\\
&v\in L_{loc}^{\infty}(\overline\Om\times[0,\infty))\cap L_{loc}^{1}([0,\infty);W^{1,1}(\Om))
\end{split}
\right.
\eeno
satisfying
\beno
u^m\na v\in L_{loc}^{1}\left(\overline\Om\times[0,\infty);\R^n\right) \quad{\rm{and}}\quad
u^mv\in L_{loc}^{1}\left(\overline\Om\times[0,\infty)\right)
\eeno
is a global weak solution of \eqref{s} if
\begin{align}
-\int_{0}^{\infty}\iom u\varphi_t-\iom u_0\varphi(\cdot,0)&=\f{1}{m}\int_{0}^{\infty}\iom u^m\na v\cdot\na\varphi+\f{1}{m}\int_{0}^{\infty}\iom u^mv\Delta \varphi\no\\
&\quad+\int_{0}^{\infty}\iom f(u)v\na v\cdot\na\varphi+\ell\int_{0}^{\infty}\iom uv\varphi\label{de3}
\end{align}
for all $\varphi\in C_{0}^{\infty}\left(\overline\Om\times[0,\infty)\right)$ fulfilling $\f{\partial\varphi}{\partial\nu}=0$ on $\partial\Om\times(0,\infty)$, as well as
\beq\label{de4}
\int_{0}^{\infty}\iom v\varphi_t+\iom v_0\varphi(\cdot,0)=\int_{0}^{\infty}\iom \na v\cdot\na\varphi+\int_{0}^{\infty}\iom uv\varphi
\eeq
for each $\varphi\in C_{0}^{\infty}\left(\overline\Om\times[0,\infty)\right)$.
\end{Definition}

Now we can state our main results.
\begin{Theorem}\label{global existence}
Let $\Omega\subset \R^2$ be a bounded convex domain with smooth boundary, and let $\ell\ge0$. Suppose that the initial data $u_0$ and $v_0$ satisfy \eqref{indata}. Then if one of the following cases holds:\\
(i) $1\le m<2$, $f$ fulfills \eqref{f1} with $m-1<\al< m$;\\
(ii) $2\le m<3$, $f$ fulfills \eqref{f2} with $m-1<\al< \f{m}{2}+1$,\\
there exist functions
\begin{equation}\label{reg1}
\left\{
\begin{split}
&u\in C^0(\overline\Om\times[0,\infty))
~~~~and\\
&v\in C^0(\overline\Om\times[0,\infty))\cap C^{2,1}(\overline\Om\times(0,\infty))
\end{split}
\right.
\end{equation}
such that $(u,v)$ forms a global weak solution of \eqref{s} in the sense of Definition \ref{de}. Moreover, if\\
(iii) $3\le m<4$, $f$ fulfills \eqref{f2} with $m-1<\al<\f{m}{2}+1$ and $u_0>0$ in $\overline{\Om}$,\\
there exist functions
\begin{equation}\label{reg1}
\left\{
\begin{split}
&u\in \cap_{q\ge1}C^0([0,\infty);{W^{1,q}}(\Omega)\big)\cap C^{2,1}(\bar\Omega\times (0,\infty))
~~~~and\\
&v\in \cap_{q\ge1}C^0([0,\infty);{W^{1,q}}(\Omega)\big)\cap C^{2,1}(\bar\Omega\times (0,\infty))
\end{split}
\right.
\end{equation}
such that $(u,v)$ forms a global classical solution of \eqref{s}.

Furthermore, $u\ge 0$ and $v>0$ in $\overline\Om\times[0,\infty)$, {and}
\beq\label{th2r}
\|u(\cdot,t)\|_{L^{\infty}(\Om)}+\|v(\cdot,t)\|_{W^{1,\infty}(\Om)}<\infty
\qquad for~a.e.~t>0.
\eeq
\end{Theorem}

\begin{Remark}
{\rm{Compared to the results of Theorem 1.2 in \cite{wu2024preprint}, we improve the $L^p$ regularity of $u$ to $L^{\infty}$, which essentially provides the possibility for the following result on asymptotic stability.}}
\end{Remark}

\begin{Theorem}\label{largetime}
Suppose that the assumptions in Theorem \ref{global existence} are satisfied, and let $(u,v)$ be as accordingly given by Theorem \ref{global existence}. Then there exists $u_{\infty}\in C^0(\overline\Om)$ such that
\beq\label{con}
u(\cdot,t)\rightarrow u_{\infty}~~and~~v(\cdot,t)\rightarrow 0~~in~~L^{\infty}(\Om)~~as ~t\rightarrow\infty.
\eeq
Here the limit function satisfies $u_{\infty}=w(\cdot,1)$ with $w\in C^0{(\overline\Om\times[0,1])}$ being a weak solution of
\beno
\left\{
\begin{split}
&w_{\tau}=\na\cdot\big(a(x,\tau)w^{m-1}\na w\big)-\na\cdot\big(b(x,\tau)f(w)\big)+\ell a(x,\tau)w,&&x\in\Omega,\,\tau\in(0,1),\\
&\nabla w\cdot\nu=0, &&x\in\p\Omega,\,\tau\in(0,1),\\
&w(x,0)=u_0(x), &&x\in\Omega,\\
\end{split}
\right.
\eeno
in the sense that
\bea{th-de3}
-\int_{0}^{1}\iom w\varphi_t-\iom w_0\varphi(\cdot,0)
&=&\f{1}{m}\int_{0}^{1}\iom w^m\na a(x,\tau)\cdot\na\varphi
+\f{1}{m}\int_{0}^{1}\iom a(x,\tau)w^m \Delta \varphi\no\\
&& +\int_{0}^{1}\iom b(x,\tau)f(w)\cdot\na\varphi
+\ell\int_{0}^{1}\iom a(x,\tau)w\varphi
\eea
for all $\varphi\in C_{0}^{\infty}\left(\overline\Om\times[0,1)\right)$ fulfilling $\f{\partial\varphi}{\partial\nu}=0$ on $\partial\Om\times(0,1)$, where
\beno
a(x,\tau):=L\cdot\f{v(x,t)}{\|v(\cdot,t)\|_{{L^{\infty}}(\Om)}} ~~{{and}}~~b(x,\tau):=L\cdot\f{v(x,t)\na v(x,t)}{\|v(\cdot,t)\|_{{L^{\infty}}(\Om)}}, ~~(x,\tau)\in \Omega\times(0,1)~~and~~t=\phi^{-1}(\tau)
\eeno
with
\beno
{L:=\int_0^{\infty}\|v(\cdot,s)\|_{{L^{\infty}}(\Om)}ds~~
and~~\phi(t):=\f{1}{L}\cdot\int_0^t\|v(\cdot,s)\|_{{L^{\infty}}(\Om)}ds,~~t\ge 0}
\eeno
are such that there exists $C>0$ satisfying
\beq\label{th-de4}
\f{1}{C}\le a(x,\tau)\le C ~~{and}~~|b(x,\tau)|\le C~~for~all~~(x,\tau)\in \Omega\times(0,1).
\eeq
\end{Theorem}
\begin{Remark}
{\rm{Although we only state the asymptotic behavior in the large time limit for spatially two-dimensional version of \eqref{s}, the approach developed in this paper appears to be available to \eqref{s} and its variants in higher dimensions.}}
\end{Remark}

\section{Some preliminaries}
Similar to the approximating procedure used in \cite{wu2024preprint}, we consider the regularized variant of \eqref{s} given by
\begin{equation}\label{s1}
\left\{
\begin{split}
&{\uep}_t=\nabla\cdot\left(u_{\e}^{m-1}{\vep}\nabla {\uep}\right)-\nabla\cdot\left(f({\uep}){\vep}\nabla {\vep}\right)+\ell {\uep}{\vep},&&x\in\Omega,\,t>0,\\
&{\vep}_t=\Delta {\vep}-{\uep}{\vep}, &&x\in\Omega,\,t>0,\\
&\f{\p\uep}{\p\nu}=\f{\p\vep}{\p\nu}=0,&& x\in\partial\Omega,\,t>0,\\
&\uep(x,0)=u_{0\e}(x),~~v(x,0)=v_{0\e}(x):=v_0(x),&&x\in\Omega
\end{split}
\right.
\end{equation}
with $\e\in(0,1)$, where {$u_{0\e}(x)$ depending on $m$ is defined by
\begin{equation}\label{u0e}
u_{0\e}(x):=\left\{
\begin{split}
&u_0(x)+\e, &&1\le m<3,\\
&u_0(x), &&3\le m<4.\\
\end{split}
\right.
\end{equation}}

The following lemma is a direct consequence of Lemma 2.1, Lemma 5.2, and Lemma 5.6 in \cite{wu2024preprint}.
\begin{Lemma}\label{lelocal}
Suppose that the assumptions in Theorem \ref{global existence} are satisfied. Then for each $\e\in(0,1)$, there exists at least one pair $(\uep,\vep)$ of functions
\beq\label{local}
\left\{
\begin{split}
&\uep\in \cap_{q\ge1} C^0\left([0,\infty);W^{1,q}(\Omega)\right)\cap C^{2,1}\left(\overline\Omega\times (0,\infty)\right)\\
&\vep\in \cap_{q\ge1} C^0\left([0,\infty);W^{1,q}(\Omega)\right)\cap C^{2,1}\left(\overline\Omega\times (0,\infty)\right)
\end{split}
\right.
\eeq
such that $\uep,\vep>0$ in $\overline\Omega\times (0,\infty)$, and there exist $(\e_j)_{j\in \N}\subset(0,1)$ fulfilling $\e_j\rightarrow 0$ as $j\rightarrow \infty$ and a pair $(u,v)$ of nonnegative functions
\beq\label{le21-1}
\left\{
\begin{split}
&u\in L^{\infty}((0,\infty);L^2(\Om))\qquad{\rm{and}}\\
&v\in L^{\infty}((0,\infty);W^{1,\infty}(\Om))
\end{split}
\right.
\eeq
such that $(u,v)$ forms a global weak solution of \eqref{s} in the sense of Definition \ref{de}, and that
\beq\label{uvc-1}
{u_{\e_j}\rightarrow u~~and~~ v_{\e_j}\rightarrow v
\qquad a.e.~in~\Om\times(0,\infty)
\quad as ~~\e_j\rightarrow 0.}
\eeq
\end{Lemma}
From now on, without further explicit mention, it is assumed that $u_0$ and $v_0$ always fulfill \eqref{indata}. Now we further introduce some elementary boundedness properties of the approximate system.
\begin{Lemma}\label{elebdd}
Suppose that $(\uep,\vep)$ is {given} by Lemma \ref{lelocal}. Then for any $\e\in(0,1)$, we have
\beq\label{vin}
\|\vep(\cdot,t)\|_{L^{\infty}({\Om})}\leq \|\vep(\cdot,t_0)\|_{L^{\infty}({\Om})}~~for~all~t_0\ge0~and~t>t_0
\eeq
and
\beq\label{u1}
\iom u_{0\e}\le\iom \uep(\cdot,t)\le\iom u_{0\e}+\ell\iom v_{0\e}~~for~all~t>0
\eeq
as well as
\beq\label{uv1}
\int_{t_0}^{\infty}\iom \uep\vep\leq \iom \vep(\cdot,t_0)~~~~for~all~t_0\ge0.
\eeq
Moreover, for any $p\ge1$ {and $\e\in(0,1)$}, there exists $C>0$ such that
\beq\label{nav6v5}
\int_{0}^{\infty}\iom \f{|\na \vep|^{6}}{\vep^{5}}\leq C
\eeq
and
\beq\label{up}
\iom \uep^p(\cdot,t)\leq C
\eeq
as well as
\beq\label{nav}
\|\na\vep(\cdot,t)\|_{L^{\infty}(\Om)}\le C
~~for~all~t>0.
\eeq
\end{Lemma}
\Bp The properties in \eqref{vin}-\eqref{uv1} can be obtained from Lemma 2.1 and Lemma 5.2 in \cite{wu2024preprint}. As an implication of Lemmata 4.5-4.7 in \cite{wu2024preprint}, we can claim that there exists $c_1>0$ such that
\beno
\int_0^{\infty}\iom \f{|\na \vep|^{2}}{v_{\e}}|D^2\ln\vep|^2\le c_1 \qmf~\e\in(0,1),
\eeno
which in conjunction with \cite[Lemma 2.2]{wu2024preprint} implies \eqref{nav6v5}. Finally, \eqref{up} and \eqref{nav} are direct results of Lemma 4.9 and Lemma 5.1, respectively, in \cite{wu2024preprint}.
\Ep

\section{Uniform $L^{\infty}$ boundedness of $u$ and the proof of Theorem \ref{global existence}}

The following inequality is taken from \cite[Lemma 6.2]{winkler2024jde}, which plays a critical role in the iterative argument leading to the $L^{\infty}$ bound for $\uep$. For completeness, we include the proof here.

\begin{Lemma}\label{lecru}
Let $\Om\subset\R^2$ and $p_*>2$. Then there exist $\kappa=\kappa(p_*)>0$ and $K=K(p_*)>0$ such that for any choice of $p\ge p_*$ and $\eta\in(0,1]$,
\beq\label{cru}
\iom \varphi^{p+1}\psi\leq \eta\iom\varphi^{p-1}\psi|\na \varphi|^2+\eta\cdot\left\{\iom\varphi^{\f{p}{2}}\right\}^{\f{2(p+1)}{p}}\cdot\iom\f{|\na \psi|^6}{\psi^5} +K\eta^{-\kappa}p^{2\kappa}\cdot\left\{\iom\varphi^{\f{p}{2}}\right\}^{2}\cdot\iom \varphi\psi
\eeq
is valid for arbitrary positive functions $\varphi\in C^1(\overline\Om)$ and $\psi\in C^1(\overline\Om)$.
\end{Lemma}
\Bp
As $p_*>2$, we have
\begin{equation}\label{ch2-n1}
q\equiv q(p_*):=\f{6p_*}{5p_*+2}>1,
\end{equation}
so that the Gagliardo-Nirenberg inequality in the two-dimensional domain $\Om$ provides $c_1>0$ such that
\begin{equation*}
\big\|  \rho \big\|_{L^2(\Om)}^2
\le
c_1\big\| \na\rho \big\|_{L^{q}(\Om)}^{\f{2q}{2q-1}}
\big\| \rho \big\|_{L^{\f{2}{3}}(\Om)}^{\f{2q-2}{2q-1}}
+{c_1}\big\| \rho \big\|_{L^{\f{2}{3}}(\Om)}^{2}
\qmf \rho\in C^1(\overline\Om).
\end{equation*}
Given $p\ge p_*$, $\eta\in(0,1]$ as well as $0<\varphi\in C^1(\overline\Om)$ and $0<\psi\in C^1(\overline\Om)$, we thus obtain that 
\bea{ch2-n2}
\iom \varphi^{p+1}\psi 
&=& \big\| \varphi^{\f{p+1}{2}}\psi^{\f{1}{2}} \big\|_{L^2(\Om)}^2\no\\
&\le&
c_1\big\|\na \left(\varphi^{\f{p+1}{2}}\psi^{\f{1}{2}} \right)\big\|_{L^{q}(\Om)}^{\f{2q}{2q-1}}
\big\| \varphi^{\f{p+1}{2}}\psi^{\f{1}{2}} \big\|_{L^{\f{2}{3}}(\Om)}^{\f{2q-2}{2q-1}}
+{c_1} \big\| \varphi^{\f{p+1}{2}}\psi^{\f{1}{2}} \big\|_{L^{\f{2}{3}}(\Om)}^{2},
\eea
where writing 
\begin{equation}\label{ch2-n3}
c_2:=\max\left\{1, ~|\Om|^{\f{2}{q}} \right\}
\quad {\rm{and}}\quad
\delta\equiv\delta(p,\eta):
=\min\bigg\{ \f{\eta} {(p+1)^2 |\Om|^{\f{2-q}{q}}},~\f{\eta^{\f{1}{3}}}{c_2} \bigg\},
\end{equation}
using Young's inequality we find that
\bea{ch2-n4}
&&\hspace{-15mm}
c_1 \left\|\na \left(\varphi^{\f{p+1}{2}}\psi^{\f{1}{2}} \right)\right\|_{L^{q}(\Om)}^{\f{2q}{2q-1}}\big\| \varphi^{\f{p+1}{2}}\psi^{\f{1}{2}} \big\|_{L^{\f{2}{3}}(\Om)}^{\f{2q-2}{2q-1}}\no\\
&=& \left\{ \delta \left\| \na \left( \varphi^{\f{p+1}{2}} \psi^{\f{1}{2}} \right) \right\|_{L^{q}(\Om)}^2 \right\}^{\f{q}{2q-1}}
\cdot c_1\delta^{-\f{q}{2q-1}} \big\| \varphi^{\f{p+1}{2}}\psi^{\f{1}{2}} \big\|_{L^{\f{2}{3}}(\Om)}^{\f{2q-2}{2q-1}}\no\\
&\le& \delta \big\|\na \left(\varphi^{\f{p+1}{2}}\psi^{\f{1}{2}} \right)\big\|_{L^{q}(\Om)}^2
+c_1^{\f{2q-1}{q-1}}\delta^{-\f{q}{q-1}} \big\| \varphi^{\f{p+1}{2}}\psi^{\f{1}{2}} \big\|_{L^{\f{2}{3}}(\Om)}^{2}.
\eea
Here, once more by Young's inequality, 
\bea{ch2-n5}
\delta \left \| \na \left(\varphi^{\f{p+1}{2}}\psi^{\f{1}{2}} \right) \right\|_{L^{q}(\Om)}^2
&=& \delta \left \| \f{p+1}{2} \varphi^{\f{p-1}{2}} \psi^{\f{1}{2}} \nabla\varphi 
              + \f{1}{2} \varphi^{\f{p+1}{2}} \psi^{-\f{1}{2}} \nabla \psi 
   \right \|_{L^q(\Om)}^2\no\\
&\le& \delta \cdot \left\{ 
     \f{p+1}{2} \left\| \varphi^{\f{p-1}{2}}\psi^{\f{1}{2}} \nabla\varphi \right\|_{L^{q}(\Om)}
     + \f{1}{2} \left\|\varphi^{\f{p+1}{2}}\psi^{-\f{1}{2}}\nabla\psi \right\|_{L^q(\Om)}
                     \right\}^2\no\\
&\le& \f{(p+1)^2\delta}{2}  
\left\| \varphi^{\f{p-1}{2}}\psi^{\f{1}{2}}\nabla\varphi \right\|_{L^{q}(\Om)}^2 
+ \f{\delta}{2} \left\| \varphi^{\f{p+1}{2}}\psi^{-\f{1}{2}}\nabla\psi \right\|_{L^q(\Om)}^2,
\eea
and observing that $q<2$ we may reply on the H\"{o}lder inequality to estimate 
\bea{ch2-n6}
\f{(p+1)^2\delta}{2}  \left\| \varphi^{\f{p-1}{2}}\psi^{\f{1}{2}}\nabla\varphi \right\|_{L^{q}(\Om)}^2
&\le&  \f{(p+1)^2\delta}{2} \cdot |\Om|^{\f{2-q}{q}} 
\left\| \varphi^{\f{p-1}{2}}\psi^{\f{1}{2}}\nabla\varphi \right\|_{L^{2}(\Om)}^2\no\\
&\le& \f{\eta}{2} \iom \varphi^{p-1}\psi|\na \varphi|^2
\eea
according to the first restriction on $\delta$ contained in \eqref{ch2-n3}. Apart from that, again by means of the H\"{o}lder inequality we see that
\bea{ch2-n7}
\f{\delta}{2} \left\| \varphi^{\f{p+1}{2}}\psi^{-\f{1}{2}}\nabla\psi \right\|_{L^q(\Om)}^2
&=& \f{\delta}{2} \cdot \left\{ \iom \varphi^{\f{(p+1)q}{2}}\psi^{-\f{q}{2}}|\nabla \psi|^q
\right \}^{\f{2}{q}}\no\\
&=& \f{\delta}{2} \cdot 
\left\{ 
\iom \left( \f{|\na \psi|^6}{\psi^5} \right)^{\f{q}{6}} 
\cdot \varphi^{\f{(p+1)q}{2}}\psi^{\f{q}{3}} 
\right\}^{\f{2}{q}}\no\\
&\le & \f{\delta}{2} \cdot \left\{ \iom \f{|\na \psi|^6}{\psi^5} \right\}^{\f{1}{3}}
\cdot \left\{ \iom \varphi^{\f{3(p+1)q}{6-q}}\psi^{\f{2q}{6-q}} \right\}^{\f{6-q}{3q}}
\eea
and that here 
\bea{ch2-n8}
\left\{ \iom \varphi^{\f{3(p+1)q}{6-q}}\psi^{\f{2q}{6-q}} \right\}^{\f{6-q}{3q}}
&=& \left\{ \iom (\varphi^{p+1}\psi)^{\f{2q}{6-q}}
\cdot \psi^{\f{(p+1)q}{6-q}} \right\}^{\f{6-q}{3q}}\no\\
&\le&  \left\{ \iom \varphi^{p+1} \psi \right\}^{\f{2}{3}}
\cdot  \left\{ \iom \varphi^{\f{(p+1)q}{6-3q}} \right\}^{\f{2-q}{q}}.
\eea
Now our definition \eqref{ch2-n1} applies in its full strength so as to assert, namely, that since $\f{d}{d\xi} \f{6\xi}{5\xi+2}\ge 0$ for all $\xi>0$, the inequality $p\ge p_*$ ensures that $q\le 
\f{6p}{5p+2}$ and hence $\f{(p+1)q}{6-3q}=\f{p+1}{\f{6}{q}-3}\le \f{p+1}{\f{5p+2}{p}-3}=\f{p}{2}$, so that a final application of the H\"older inequality shows that
\bea{ch2-n9}
 \left\{ \iom \varphi^{\f{(p+1)q}{6-3q}} \right\}^{\f{2-q}{q}}
\le |\Om|^{\f{6p-5pq-2q}{3pq}} \cdot \left\{ \iom \varphi^{\f{p}{2}} \right\}^{\f{2(p+1)}{3p}}
\le  c_2 \cdot \left\{ \iom \varphi^{\f{p}{2}} \right\}^{\f{2(p+1)}{3p}}
\eea
with $c_2$ as in \eqref{ch2-n3}, because clearly $0\le \f{6p-5pq-2q}{3pq} \le \f{2}{q}$.\\
From \eqref{ch2-n7}, \eqref{ch2-n8} and \eqref{ch2-n9} we now obtain, employing Young's inequality once again, that 
\bas
\f{\delta}{2} \left\| \varphi^{\f{p+1}{2}} \psi^{-\f{1}{2}} \na \psi \right\|_{L^q(\Om)}^2
&\le& \f{c_2\delta}{2} 
\cdot \left\{ \iom \f{|\na\psi|^6}{\psi^5} \right\}^{\f{1}{3}}
\cdot \left\{ \iom \varphi^{p+1} \psi \right\}^{\f{2}{3}}
\cdot \left\{ \iom \varphi^{\f{p}{2}} \right\}^{\f{2(p+1)}{3p}}\\
&=& \left\{ \f{1}{2} \iom \varphi^{p+1} \psi \right\}^{\f{2}{3}}
\cdot \f{c_2\delta}{2^{\f{1}{3}}}
\cdot \left\{ \iom \f{|\na\psi|^6}{\psi^5} \right\}^{\f{1}{3}}
\cdot \left\{ \iom \varphi^{\f{p}{2}} \right\}^{\f{2(p+1)}{3p}}\\
&\le& \f{1}{2} \iom \varphi^{p+1}\psi
+\f{c_2^3 \delta^3}{2} 
\cdot \left\{ \iom \varphi^{\f{p}{2}} \right\}^{\f{2(p+1)}{p}}
\cdot  \iom \f{|\na\psi|^6}{\psi^5} ,
\eas
whence collecting \eqref{ch2-n2}, \eqref{ch2-n4}, \eqref{ch2-n5} and \eqref{ch2-n6} we see that
\bea{ch2-n10}
 \iom \varphi^{p+1}\psi
&\le& \f{\eta}{2}\iom \varphi^{p-1} \psi |\na \varphi|^2
+\f{1}{2} \iom \varphi^{p+1}\psi 
+\f{c_2^3 \delta^3}{2}
\cdot \left\{ \iom \varphi^{\f{p}{2}} \right\}^{\f{2(p+1)}{p}}
\cdot \iom \f{|\na\psi|^6}{\psi^5} \no\\
&& + \left( c_1^{\f{2q-1}{q-1}} \delta^{-\f{q}{q-1}}+c_1 \right)
\cdot \left\| \varphi^{\f{p+1}{2}} \psi^{\f{1}{2}} \right\|_{L^{\f{2}{3}}(\Om)}^2.
\eea
As 
\bas
\left\| \varphi^{\f{p+1}{2}} \psi^{\f{1}{2}} \right\|_{L^{\f{2}{3}}(\Om)}^2
= \left\{ \iom \varphi^{\f{p+1}{3}}\psi^{\f{1}{3}} \right\}^3
=  \left \{ \iom (\varphi\psi)^{\f{1}{3}} \cdot \varphi^{\f{p}{3}} \right\}^3
\le \left\{\iom \varphi^{\f{p}{2}} \right\}^2 \cdot \iom \varphi\psi,
\eas
this entails that
\bas
\iom \varphi^{p+1} \psi 
&\le& \eta \iom \varphi^{p-1}\psi|\na \varphi|^2 
+ c_2^3 \delta^3 
\cdot \left\{\iom \varphi^{\f{p}{2}} \right\}^{\f{2(p+1)}{p}}
\cdot \iom \f{|\na\psi|^6}{\psi^5} \\
&& + 2\left( c_1^{\f{2q-1}{q-1}} \delta^{-\f{q}{q-1}}+c_1 \right)
\cdot \left\{ \iom \varphi^{\f{p}{2}} \right\}^2 \cdot \iom \varphi\psi ,
\eas
and thereby establishes \eqref{ch2-n10} with 
\bas
\kappa\equiv\kappa(p_*):=\f{q}{q-1}
\quad {\rm{and}} \quad
K\equiv K(p_*):= 2c_1^{\f{2q-1}{q-1}} \cdot \max\left\{ \left(4|\Om|^{\f{2-q}{q}}\right)^{\f{q}{q-1}},~c_2^{\f{q}{q-1}} \right\} +2c_1
\eas
because the inequalities $p\ge 1$ and $\eta\le 1$ warrant that, by \eqref{ch2-n3},
\bas
2 c_1^{\f{2q-1}{q-1}}\delta^{-\f{q}{q-1}}
&=& 2 c_1^{\f{2q-1}{q-1}} 
\cdot \max \left\{ \left( \f{(p+1)^2|\Om|^{\f{2-q}{q}}}{\eta} \right)^{\f{q}{q-1}}
,~\left( \f{c_2}{\eta^{\f{1}{3}}} \right)^{\f{q}{q-1}} \right\}\\
&\le& 2c_1^{\f{2q-1}{q-1}} 
\cdot \max \left\{ \left( \f{(2p)^2|\Om|^{\f{2-q}{q}}}{\eta} \right)^{\f{q}{q-1}}
,~\left( \f{c_2p^2}{\eta} \right)^{\f{q}{q-1}} \right\}\\
&\le& 2c_1^{\f{2q-1}{q-1}} \eta^{-\f{q}{q-1}} p^{\f{2q}{q-1}} 
\cdot \max \left\{ \left(4|\Om|^{\f{2-q}{q}}\right)^{\f{q}{q-1}}, ~c_2^{\f{q}{q-1}} \right\},
\eas
and that $2c_1\le 2c_1\eta^{-\f{q}{q-1}}p^{\f{2q}{q-1}}$.
\Ep

We are now able to establish the $L^{\infty}$ bound for $\uep$ by using the Moser iterative technique.
\begin{Lemma}\label{leuin}
Suppose that the assumptions in Theorem \ref{global existence} are satisfied. Let $(\uep,\vep)$ be as yielded by Lemma \ref{lelocal}. Then there exists $C>0$ such that
\beno
\|\uep(\cdot,t)\|_{L^{\infty}(\Om)}\le C\quad for~all~ t>0~and ~\e\in(0,1)
\eeno
\end{Lemma}
\Bp Take $p_0=4$, and recursively define
\beq\label{pk}
p_k:=2p_{k-1}+2-m,~~~~{k\in \{1,2,3,...\}}.
\eeq
Then it is obvious that $(p_k)_{k\in\N}$ increases and
\beq\label{pk1}
c_1\cdot 2^k\le p_k\le c_2\cdot 2^k\qmf k\in\N
\eeq
with $c_1:=p_0-(2-m)_{-}$ and $c_2:=p_0+(2-m)_{+}$. Setting
\beq\label{mk}
M_{k,\e}(T):=1+\sup_{t\in(0,T)}\iom \uep^{p_k}(\cdot,t),~~T\in(0,\infty),~~k\in\N~~\rm{and}~~\e\in(0,1),
\eeq
then {we see that} each $M_{k,\e}(T)$ is finite and we can use \eqref{up} to see the existence of $c_3>0$ independent of $T$ satisfying
\beq\label{m0}
M_{0,\e}\le c_3 \qmf \e\in(0,1).
\eeq
Now we try to estimate $M_{k,\e}(T)$ for $T\in(0,\infty)$, $k\ge1$ and $\e\in(0,1)$. By \eqref{nav}, we have
\beno
\|\na\vep(\cdot,t)\|_{L^{\infty}(\Om)}\le c_4 \qmf t>0~ {{\rm{and}}~\e\in(0,1)}.
\eeno
with $c_4>0$. And recalling \eqref{f1} and \eqref{f2}, we can claim that there exists $c_5>0$ such that
\beno
f^2(\uep)\le c_5C_f^2\left(\uep^{2\al}+\uep^2\right) \qmf \e\in(0,1).
\eeno
Thus, {after} testing the first equation in \eqref{s1} by $p_k\uep^{p_k-1}$ and integrating by parts, the Young inequality along with the boundary conditions entails that
\bea{upk1}
\f{d}{dt}\iom \uep^{p_k}
&=& {p_k}\iom \uep^{{p_k}-1}\nabla\cdot\left(u_{\e}^{m-1}{\vep}\nabla {\uep}\right)-{p_k}\iom \uep^{{p_k}-1}\nabla\cdot\big(f({\uep}){\vep}\nabla {\vep}\big)+{p_k}\ell\iom \uep^{{p_k}}\vep \no\\
&=& -{p_k}({p_k}-1)\iom \uep^{{p_k}+m-3}\vep|\na \uep|^2+{p_k}({p_k}-1)\iom \uep^{{p_k}-2}f(\uep)\vep\na \uep\cdot\na\vep+{p_k}\ell\iom \uep^{{p_k}}\vep\no\\
&\le& -\f{{p_k}({p_k}-1)}{2} \iom \uep^{{p_k}+m-3}\vep|\na \uep|^2+\f{{p_k}({p_k}-1)}{2}\iom \uep^{{p_k}-m-1}f^2(\uep)\vep|\na \vep|^2+{p_k}\ell\iom \uep^{{p_k}}\vep\no\\
&\le& -\f{p_k^2}{4} \iom \uep^{{p_k}+m-3}\vep|\na \uep|^2
+c_4^2c_5C_f^2{p_k^2}\left\{\iom \uep^{{p_k}-m+1}\vep
+\iom \uep^{{p_k}+2\al-m-1}\vep\right\}
+p_k\ell \iom \uep^{{p_k}}\vep\no\\
&\le& -\f{p_k^2}{4} \iom \uep^{{p_k}+m-3}\vep|\na \uep|^2
+(2c_4^2c_5C_f^2+\ell)p_k^2\iom \uep^{{p_k}+m-1}\vep\no\\
&& +(2c_4^2c_5C_f^2+\ell)p_k^2 \iom \uep\vep
\qmf t>0~ \rm{and}~\e\in(0,1),
\eea
where we also use the facts that $1<{p_k}+2\al-m-1<{p_k}+m-1$ and $1<{p_k}-m+1\le p_k\le{p_k}+m-1$ guaranteed by our restrictions $1\le m<4$ and $m-1<\al<m$.

\quad Due to $p_k\ge 4$ for all $k\ge 1$, we have ${p_k}+m-1>4$, and thus Lemma \ref{lecru} with taking $p_*:=3$ infers the existence of $\kappa>0$ and $K>0$ satisfying
\bea{upk2}
&& \hspace{-15mm}
(2c_4^2c_5C_f^2+\ell)p_k^2\iom \uep^{{p_k}+m-1}\vep\no\\
&\le& \f{p_k^2}{4}\iom\uep^{{p_k}+m-3}\vep|\na \uep|^2
+\f{p_k^2}{4}\cdot\left\{\iom\uep^
{\f{{p_k}+m-2}{2}}\right\}^{\f{2({p_k}+m-1)}
{{p_k+m-2}}}\cdot\iom\f{|\na \vep|^6}{\vep^5}\no\\ && +4^{\kappa}(2c_4^2c_5C_f^2+\ell)^{\kappa+1}Kp_k^2(p_k+m-2)^{2\kappa}\cdot
\left\{\iom\uep^{\f{{p_k}+m-2}{2}}\right\}^{2}
\cdot\iom \uep\vep
\eea
for all $t>0$ and $\e\in(0,1)$. Combining \eqref{upk1} and \eqref{upk2}, we have
\bea{upk3}
\f{d}{dt}\iom \uep^{p_k}
&\le& \f{p_k^2}{4}\cdot\left\{\iom\uep^
{\f{{p_k}+m-2}{2}}\right\}^{\f{2({p_k}+m-1)}
{{p_k+m-2}}}\cdot\iom\f{|\na \vep|^6}{\vep^5}\no\\
&& +c_7p_k^2(p_k+m-2)^{2\kappa}\cdot
\left\{\iom\uep^{\f{{p_k}+m-2}{2}}\right\}^{2}
\cdot\iom \uep\vep
+c_6p_k^2 \iom \uep\vep
\eea
for all $t>0$ and $\e\in(0,1)$ with $c_6:=2c_4^2c_5C_f^2+\ell$ and $c_7:=4^{\kappa}c_6^{\kappa+1}K$. From \eqref{pk} and \eqref{pk1}, we have
\beno
c_1\cdot 2^k\le{p_k}+m-2=2p_{k-1}\le c_2\cdot 2^k {\qmf k\in \{1,2,3,...\},}
\eeno
which together with \eqref{pk1}, \eqref{mk} and \eqref{upk3} implies that for all $t\in(0,T)$, $T\in(0,\infty)$ and $\e\in (0,1)$, we have
\bas
\f{d}{dt}\iom \uep^{p_k}
&\le& \f{(2^kc_2)^2}{4}\cdot\left\{\iom\uep^{p_{k-1}}\right\}^{2+\f{2}{p_k+m-2}}
\cdot\iom\f{|\na \vep|^6}{\vep^5}\no\\
&& +c_7(2^kc_2)^{2\kappa+2}\cdot
\left\{\iom\uep^{p_{k-1}}\right\}^{2}
\cdot\iom \uep\vep
+c_6(2^kc_2)^{2}\iom \uep\vep\no\\
&\le&  \f{c_2^2}{4} (2^k)^{2} M_{k-1,\e}^{2+\f{2}{c_1}\cdot 2^{-k}}(T)
\cdot\iom\f{|\na \vep|^6}{\vep^5}
+c_7c_2^{2\kappa+2}(2^k)^{2\kappa+2} M_{k-1,\e}^{2}(T)
\cdot\iom \uep\vep\no\\
&& +c_6c_2^{2}(2^k)^{2}\iom \uep\vep\no\\
&\le& c_8(2^k)^{2\kappa+2} M_{k-1,\e}^{2+\f{2}{c_1}\cdot 2^{-k}}(T)
\cdot\left\{\iom\f{|\na \vep|^6}{\vep^5}+\iom \uep\vep\right\}
\qmf t>0~ \rm{and}~\e\in(0,1)
\eas
with $c_8:=\f{c_2^2}{4}+c_7c_2^{2\kappa+2}+c_6c_2^2$. Integrating this in time, we see that for all $t\in(0,T)$, $T\in (0,\infty)$ and $\e\in(0,1)$,
\begin{align}\label{upk4}
\iom \uep^{p_k}\le c_8c_9(2^k)^{2\kappa+2} M_{k-1,\e}^{2+\f{2}{c_1}\cdot 2^{-k}}(T)+\iom (u_0+1)^{p_k},
\end{align}
where
\beno
c_9:={\sup_{\e\in(0,1)}}\left\{\int_0^{\infty}\iom\f{|\na \vep|^6}{\vep^5}
+\int_0^{\infty}\iom \uep\vep\right\}<\infty
\eeno
warranted by \eqref{uv1} and \eqref{nav6v5}. If we write
\beno
a:=1+\left(|\Om|+1\right)\|u_0+1\|_{L^{\infty}(\Om)}^{c_2}
~~\mathrm{and}~~b:=(c_8c_9+1)\cdot 2^{2\kappa+2},
\eeno
then it is easy to verify that
\beno
\|u_0+1\|_{L^{\infty}(\Om)}^{p_k}\cdot|\Om|+1\le 1+\left(|\Om|+1\right)^{2^k}\|u_0+1\|_{L^{\infty}(\Om)}^{c_2\cdot2^k}\le a^{2^k}
\eeno
and
\beno
 c_8c_9(2^k)^{2\kappa+2} \le (c_8c_9+1)^k(2^{2\kappa+2})^k= b^k.
\eeno
That is, we can further conclude from \eqref{upk4} and \eqref{mk} that
\bas
M_{k,\e}(T)
&\le& c_8c_9(2^k)^{2\kappa+2} M_{k-1,\e}^{2+\f{2}{c_1}\cdot 2^{-k}}(T)+ \|u_0+1\|_{L^{\infty}(\Om)}^{p_k}\cdot|\Om|+1\no\\
&\le& b^k M_{k-1,\e}^{2+\f{2}{c_1}\cdot 2^{-k}}(T)+a^{2^k}.
\eas
Since $k\ge 1$ is arbitrary here, together with \eqref{pk1} and \eqref{m0}, we may use \cite[Lemma 6.3]{winkler2024jde} to claim that
\bas
\|\uep(\cdot,t)\|_{L^{\infty}(\Om)}^{c_1}
&=&\liminf_{k\rightarrow\infty}\left\{\iom \uep^{p_k}(\cdot,t)\right\}^{\f{c_1}{p_k}}
\le \liminf_{k\rightarrow\infty}M_{k,\e}^{\f{c_1}{p_k}}(T)\no\\
&\le& \liminf_{k\rightarrow\infty}M_{k,\e}^{\f{1}{2^k}}(T)
\le (2\sqrt{2}b^3a^{1+\f{1}{c_1}}c_3)^{e^{\f{1}{c_1}}}
\eas
for all $t\in(0,T),~T\in(0,\infty)$ and $\e\in(0,1)$. This clearly proves the lemma.
\Ep

We finally show that $\uep$ and $\vep$ enjoy higher regularities.
\begin{Lemma}\label{uvh}
Suppose that the assumptions in Theorem \ref{global existence} are satisfied. Let $(\uep,\vep)$ be as yielded by Lemma \ref{lelocal}. Then for any $T_1>0$, there exist $\theta_1=\theta_1(T_1)\in(0,1)$ {and $C_1(T_1)>0$} such that
\beq\label{uholder}
{\|\uep\|_{C^{\theta_1,\f{\theta_1}{2}}\left(\overline\Om\times[0,T_1]\right)}
\le C_1(T_1) ~~~~for~all~\e\in(0,1)}
\eeq
and
 \beq\label{vholder}
{\|\vep\|_{C^{\theta_1,\f{\theta_1}{2}}\left(\overline\Om\times[0,T_1]\right)}
\le C_1(T_1) ~~~~for~all~\e\in(0,1)}.
\eeq
In addition, for each $\tau>0$ and any $T_2>\tau$, there exist $\theta_2=\theta_2(\tau, T_2)\in(0,1)$ {and $C_2(\tau, T_2)>0$} such that
\beq\label{vschauder}
{\|\vep\|_{C^{2+\theta_2,1+\f{\theta_2}{2}}\left(\overline\Om\times[\tau,T_2]\right)}
 \le  C_2(\tau, T_2) ~~~~for~all~\e\in(0,1)}.
\eeq
\end{Lemma}

\Bp
From Lemma \ref{leuin}, we see the existence of positive constant $c_1$ such that
\beno
\uep(x,t)\le c_1 \qmf x\in \Om,~t>0~\mathrm{and}~\e\in(0,1),
\eeno
which yields from the second equation in \eqref{s1} that for all $\e\in(0,1)$,
\beno
v_{\e t}\ge \Delta \vep-c_1\vep \qquad \mathrm{in}~\Om\times (0,\infty).
\eeno
Then by the comparison principle we get that for all $\e\in(0,1)$,
\beq\label{vlower}
\vep(x,t)\ge c_2e^{-c_1t}  \qmf x\in \Om,~t>0~\mathrm{and}~\e\in(0,1)
\eeq
with $c_2=\inf_{\Om}v_0>0$ due to the strict positivity of $v_0$ asserted by \eqref{indata}.

Now we rewrite the first equation in \eqref{s1} in the following form
\beno
{\uep}_t=\nabla\cdot A_{\e}(x,t,\uep,\na\uep)+B_{\e}(x,t,\uep),\qquad x\in\Omega,\,t>0
\eeno
with
\begin{align*}
&A_{\e}(x,t,\uep,\na\uep):=\vep(x,t)\uep^{m-1}(x,t)\na\uep(x,t)-f(\uep(x,t))\vep(x,t)
\na\vep(x,t)\qquad {\rm{and}}\no\\
&B_{\e}(x,t)=\ell\uep(x,t)\vep(x,t),\qquad\qquad\qquad\qquad (x,t)\in \Om\times(0,\infty).
\end{align*}
Recalling \eqref{f1} and \eqref{f2}, when $1\le m<2$, the assumption $m-1<\al<m$ leads
\beno
\f{f^2(\uep)}{u^{m-1}_{\e}}\le C^2_f \uep^{3-m} (\uep+1)^{2\al-2}
\le C^2_f (\uep+1)^{2\al-m+1}\le C^2_f(c_1+1)^{3} \qmf \e\in(0,1),
\eeno
and when $2\le m<4$, the assumption $m-1<\al<\f{m}{2}+1$ results
\beno
\f{f^2(\uep)}{u^{m-1}_{\e}}\le C^2_f \uep^{2\al-m+1}
\le C^2_f \uep^3\le C^2_fc_1^{3}  \qmf \e\in(0,1).
\eeno
Thus, the Young inequality combined with \eqref{vlower}, \eqref{vin} and \eqref{nav} yields that for each $T>0$, there exists $c_3>0$ such that
\bas
A_{\e}(x,t,\uep,\na\uep)\cdot\na\uep
&=& \vep\uep^{m-1}|\na\uep|^2
-f(\uep)\vep\na\vep\cdot\na\uep\no\\
&\ge& \f{1}{2}\vep\uep^{m-1}|\na\uep|^2
-\f{f^2(\uep)}{u^{m-1}_{\e}}\vep|\na\vep|^2\no\\
&\ge& \f{c_2}{2}\cdot e^{-c_1T}\uep^{m-1}|\na\uep|^2
-c_3 \qmf(x,t)\in \Om\times(0,T).
\eas
as well as
\begin{align*}
&|A_{\e}(x,t,\uep,\na\uep)|\le c_3\uep^{m-1}|\na\uep|
+c_3u^{\f{m-1}{2}}_{\e},\qquad {\rm{and}}\no\\
&|B_{\e}(x,t)|\le \ell c_1c_3 \qquad \qmf(x,t)\in \Om\times(0,T).
\end{align*}
We may invoke the H\"{o}lder estimates in \cite{holder1993} to obtain \eqref{uholder}.
{The property in} \eqref{vholder} can be achieved by proceeding with a similar but simpler argument on the second equation in \eqref{s}. Furthermore, the parabolic Schauder theory in \cite{LSU1968} is applicable to get \eqref{vschauder}.
\Ep

\noindent{\bf Proof of Theorem \ref{global existence}}. Let $(u, v)$ and $(\e_j)_{j\in \N}$ be as in Lemma \ref{lelocal}. Using the Arzel\`{a}-Ascoli theorem together with \eqref{uholder}-\eqref{vschauder}, as $\e=\e_j\rightarrow 0$ we have
\beq\label{uc}
 \uep\rightarrow u  \quad \mathrm{in}\quad C^0_{loc}\left(\overline\Om\times[0,\infty)\right)
\eeq
and
\beq\label{vc}
 \vep\rightarrow v  \quad \mathrm{in}\quad C^0_{loc}
 \left(\overline\Om\times[0,\infty)\right)\cap C^{2,1}_{loc}\left(\overline\Om\times(0,\infty)\right),
\eeq
which in conjunction with Lemma \ref{leuin} reveals that
\beno
\left\{
\begin{split}
&u\in C^0(\overline\Om\times[0,\infty))\cap L^{\infty}(\Om\times(0,\infty))~~~~\mathrm{and}\\
&v\in C^0(\overline\Om\times[0,\infty))\cap C^{2,1}(\overline\Om\times(0,\infty)).
\end{split}
\right.
\eeno
This together with \eqref{local} and \eqref{le21-1} completes the proof. $\hfill\blacksquare$

\section{Harnack-type inequality and the proof of Theorem \ref{largetime}}

Based on the $L^{\infty}$ boundedness of the first component, we can proceed to derive a Harnack-type inequality for the second component $v$, which is of essential importance to the subsequent outcome regarding asymptotic stability.
\begin{Lemma}\label{harnack}
Suppose that the assumptions in Theorem \ref{global existence} are satisfied. Let $(\uep,\vep)$ be as yielded by Lemma \ref{lelocal}. Then there exists $\lam>0$ such that
\beq\label{eq-harnack}
\vep(x,t)\ge \lam \|\vep(\cdot,t)\|_{L^{\infty}(\Om)} \quad for~all~x\in\Om,~{t>0}~and~\e\in(0,1).
\eeq
\end{Lemma}

\Bp
According to Lemma \ref{leuin}, there exists $c_1>0$ such that
\beno
{\|\uep\|_{L^{\infty}(\Om\times(0,\infty))}
\le c_1~~\qmf \e\in(0,1)},
\eeno
which together with the second equation in \eqref{s1} makes \cite[Lemma 2.5 ]{huska2006jde} become applicable so as to deduce with $\lam_*>0$ we have
\beno
\vep(x,t)\ge \lam_*\|\vep(\cdot,t)\|_{L^{\infty}(\Om)} \qmf x\in\Om,~{t> 1}~{\rm{and}}~\e\in(0,1).
\eeno
When $0<t\le 1$, \eqref{vlower} provides positive constants $c_2$ and $c_3$ such that
\beno
\vep(x,t)\ge c_2e^{-c_3t}\ge c_2e^{-c_3}  \qmf x\in \Om~\mathrm{and}~\e\in(0,1),
\eeno
whereas \eqref{vin} entails that
\beno
\|\vep(\cdot,t)\|_{L^{\infty}({\Om})}\leq \|v_0\|_{L^{\infty}({\Om})} \qmf t>0~\mathrm{and}~\e\in(0,1),
\eeno
which yields that
\beno
\vep(x,t)\ge\f{c_2e^{-c_3}}{\|v_0\|_{L^{\infty}({\Om})}}
\|\vep(\cdot,t)\|_{L^{\infty}(\Om)} \qmf x\in\Om,~{0<t\le 1}~{\rm{and}}~\e\in(0,1).
\eeno
We thereby obtain \eqref{eq-harnack} by taking $\lam=\min\left\{\lam_*, \f{c_2e^{-c_3}}{\|v_0\|_{L^{\infty}({\Om})}}\right\}$.
\Ep

With the above elliptic Harnack-type inequality at hand, we can immediately derive the following result, which is similar to that in \cite[Lemma 5.2]{li-winkler2022cpaa}.

\begin{Lemma}\label{le3-1}
{Suppose that the assumptions in Theorem \ref{global existence} are satisfied. Let $(\uep,\vep)$ be as yielded by Lemma \ref{lelocal}, and $\lam$ be taken from Lemma \ref{harnack}. Then we have}
\beq\label{eq3-1}
\int_{{0}}^{\infty}\|\vep(\cdot,t)\|_{L^{\infty}(\Om)}dt\le \f{\iom v_0}{\lam\iom u_0} \quad for~all~\e\in(0,1).
\eeq
\end{Lemma}

\Bp
Making use of \eqref{vin}-\eqref{uv1} and \eqref{eq-harnack}, for all $\e\in(0,1)$, we have
\bas
\iom v_0&\ge& \int_{{0}}^{\infty} \iom \uep\vep
\ge \lam \int_{{0}}^{\infty} \left\{\iom \uep\right\}\cdot \|\vep(\cdot,t)\|_{L^{\infty}(\Om)}ds \\
&\ge&\lam\cdot\left\{\iom u_0\right\}\cdot \int_{{0}}^{\infty} \|\vep(\cdot,t)\|_{L^{\infty}(\Om)}ds
\qmf \e\in(0,1),
\eas
which completes the proof.
\Ep

We observe that the integrability of $\|\vep(\cdot,t)\|_{L^{\infty}(\Om)}$ allows for a transformation on time scale, and the Harnack-type in Lemma \ref{harnack} will facilitate the transformed version belonging to a non-degenerate diffusion parabolic problem of porous medium type. To this end, we arrive at the following result.
\begin{Lemma}\label{le-u-infty}
With $(\uep,\vep)$ and $(\e_j)_{j\in \N}$ taken from Lemma \ref{lelocal}. Let
\beno
L_{\e}:=\int_{{0}}^{\infty}\|\vep(\cdot,t)\|_{L^{\infty}(\Om)}dt,
\qquad\e\in(\e_j)_{j\in \N},
\eeno
\beno
\tau:=\phi_{\e}(t):=\f{1}{L_{\e}}\int_{{0}}^t\|\vep(\cdot,s)\|_{{L^{\infty}}(\Om)}ds,\qquad t\ge 0
\eeno
and
\beno
w_{\e}(x,\tau):=\uep(x,\phi_{\e}^{-1}(\tau)),\qquad x\in\overline\Om,~\tau\in[0,1).
\eeno
Then we have
\beno
\left\{
\begin{split}
&w_{\e\tau}=\na\cdot\big(a_{\e}(x,\tau)w_{\e}^{m-1}\na w_{\e}\big)-\na\cdot\big(b_{\e}(x,\tau)f(w_{\e})\big)+\ell a_{\e}(x,\tau)w_{\e},&&x\in\Omega,\,\tau\in(0,1),\\
&\na w_{\e}\cdot\nu=0, &&x\in\p\Omega,\,\tau\in(0,1),\\
&w_{\e}(x,0)=u_{0}(x)+\e, &&x\in\Omega
\end{split}
\right.
\eeno
with
\beno
a_{\e}(x,\tau):=L_{\e}\cdot\f{\vep(x,t)}{\|\vep(\cdot,t)\|_{{L^{\infty}}(\Om)}} ~~{{and}}~~b_{\e}(x,\tau):=L_{\e}\cdot\f{\vep(x,t)\na v_{\e}(x,t)}{\|\vep(\cdot,t)\|_{{L^{\infty}}(\Om)}}.
\eeno
Additionally, there exists $c>0$ such that
\beq\label{eq-abe}
\f{1}{c}\le a_{\e}(x,\tau)\le c ~~{and}~~|b_{\e}(x,\tau)|\le c
~~for~all~ (x,\tau)\in \Omega\times(0,1)~~{{and}}~~\e\in(\e_j)_{j\in \N},
\eeq
and
\beq\label{eq-lecl}
L_{\e}\rightarrow L:=\int_{{0}}^{\infty}\|v(\cdot,t)\|_{{L^{\infty}}(\Om)}
~~{as}~~\e=\e_j\rightarrow 0.
\eeq
\end{Lemma}

\Bp Let $\lam$ be taken from Lemma \ref{harnack}, then \eqref{vlower} shows the existence of $c_1$ and $c_2$ fulfilling
\begin{align*}
a_{\e}(x,\tau)>\lambda L_{\e}
=\lambda\int_{{0}}^{\infty}\|\vep(\cdot,t)\|_{{L^{\infty}}(\Om)}dt
\ge \lambda c_1\int_{{0}}^{\infty}e^{-c_2 t}dt
=\f{\lambda c_1}{c_2}
\end{align*}
for all $x\in \Omega,~\tau\in(0,1)$ and $\e\in(0,1)$. The upper bounds for $a_{\e}(x,\tau)$ and $|b_{\e}(x,\tau)|$ in \eqref{eq-abe} can be resulted from Lemma \ref{le3-1} and \eqref{nav}.

It follows from \eqref{vc}, Fatou's lemma and Lemma \ref{le3-1} that,
\beq\label{eq4-1}
\int_{{0}}^{\infty}\|v(\cdot,t)\|_{L^{1}(\Om)}dt
\le|\Om|\int_{{0}}^{\infty}\|v(\cdot,t)\|_{L^{\infty}(\Om)}dt
\le \f{|\Om|\iom v_0}{\lam\iom u_0}.
\eeq
From \eqref{vc}, we know that $\|v(\cdot,t)\|_{L^{1}(\Om)}$ is uniformly continuous with respect to {$t>0$}. Therefore, an application of \cite[Lemma 3.1]{BaiWinkler2016iumj} together with \eqref{eq4-1} shows that for any $\eta>0$, there exists {$t_0>0$} such that
\beno
\|v(\cdot,t_0)\|_{L^{1}(\Om)}\le \f{\eta\lam\iom u_0}{6},
\eeno
which in conjunction with \eqref{vc} concludes that there exists $\e_*\in (0,1)$ satisfying
\beno
\|\vep(\cdot,t_0)\|_{L^{1}(\Om)}\le \f{\eta\lam\iom u_0}{3}, \qmf \e\in(\e_j)_{j\in \N}~~{\rm{with}}~~\e<\e_*,
\eeno
whence similar to the proof of Lemma \ref{le3-1}, we have
\beq\label{eq4-4}
\int_{t_0}^{\infty}\|\vep(\cdot,t)\|_{L^{\infty}(\Om)}dt
\le \f{\iom \vep(\cdot,t_0)}{\lam\iom u_0}
\le \f{\eta}{3}.
\eeq
Now we apply \eqref{vc} and Fatou's lemma once more to see that
\beq\label{eq4-5}
\int_{t_0}^{\infty}\|v(\cdot,t)\|_{L^{\infty}(\Om)}dt
\le \lim_{\e\rightarrow 0}\int_{t_0}^{\infty}\|\vep(\cdot,t)\|_{L^{\infty}(\Om)}dt
\le \f{\eta}{3}
\eeq
and moreover, we can pick $\e_{**}\in(0,1)$ fulfilling
\beq\label{eq4-6}
\|v(\cdot,t)-\vep(\cdot,t)\|_{L^{\infty}(\Om)}
\le \f{\eta}{3t_0} \qmf t\in ({{0}},t_0)~{\rm{and}}~\e\in(\e_j)_{j\in \N}~~{\rm{with}}~~\e<\e_{**}.
\eeq
Thus, combining \eqref{eq4-4}-\eqref{eq4-6} gives
\bas
|L_{\e}-L|
&=&\bigg| \int_{{0}}^{\infty}\|\vep(\cdot,t)\|_{{L^{\infty}}(\Om)}dt
-\int_{{0}}^{\infty}\|v(\cdot,t)\|_{{L^{\infty}}(\Om)}dt \bigg|\\
&\le& \int_{{0}}^{t_0}\|\vep(\cdot,t)-v(\cdot,t)\|_{{L^{\infty}}(\Om)}dt
+\int_{t_0}^{\infty}\|\vep(\cdot,t)\|_{{L^{\infty}}(\Om)}dt
+\int_{t_0}^{\infty}\|v(\cdot,t)\|_{{L^{\infty}}(\Om)}dt\\
&\le& \eta \qmf \e\in(\e_j)_{j\in \N}~~{\rm{with}}~~\e<\min\{\e_{*}, \e_{**}\}.
\eas
This thereby proves \eqref{eq-lecl}.
\Ep

\noindent{\bf Proof of Theorem \ref{largetime}}. Let $(\e_j)_{j\in \N}$ be as in Lemma \ref{lelocal}. Then according to \eqref{vc} and \eqref{eq-lecl}, we have
\beno
\phi_{\e}(t)\rightarrow \phi(t) \qmf  t>{{0}}~~\rm{as}~~\e=\e_j\rightarrow 0.
\eeno
Therefore, from \eqref{uc} and \eqref{vc} we have
\beq\label{eq-wabe}
w_{\e}(x,\tau)\rightarrow u(x,\phi^{-1}(\tau)),~ a_{\e}(x,\tau)\rightarrow a(x,\tau) ~~{\rm{and}}~~b_{\e}(x,\tau)\rightarrow b(x,\tau)
\eeq
for all $(x,\tau)\in\Om\times(0,1)$ as $\e=\e_j\rightarrow 0$.

On the other hand, thanks to $\eqref{eq-abe}$, we may {rely on} the H\"{o}lder regularity in quasilinear degenerate parabolic equations (\cite{holder1993}) to claim that there exist {$\theta\in(0,1)$ and $C>0$} such that
\beno
{\|w_{\e}\|_{C^{\theta,\f{\theta}{2}}\left(\overline\Om\times[0,1]\right)}
\le C \qmf \e\in(0,1)}
\eeno
in quite a similar manner stated in Lemma \ref{uvh}. Then by the Arzel\`{a}-Ascoli theorem, we obtain that
\beno
w_{\e}(x,\tau)\rightarrow w(x,\tau)\quad\mathrm{in}
\quad C^0\left(\overline\Om\times[0,1]\right)~~{\rm{as}}~\e=\e_j\rightarrow 0
\eeno
for some $w\in C^0\left(\overline\Om\times[0,1]\right)$. Then we can conclude that
\beno
w(x,\tau)=u(x,\phi^{-1}(\tau))\qmf (x,\tau)\in\Om\times(0,1),
\eeno
which along with the continuity of $w(\cdot,1)$ in $\overline\Om$ indicates that
\beno
u(\cdot,t)\rightarrow u_{\infty}:=w(\cdot,1)~~{\rm{in}}~~L^{\infty}(\Om)~~{\rm{as}}~t\rightarrow \infty.
\eeno
From \eqref{vin} and \eqref{vc}, we have
\begin{align*}
\|v(\cdot,t)\|_{{L^{\infty}}(\Om)}\le \|v(\cdot,t_0)\|_{{L^{\infty}}(\Om)} \qmf t_0\ge0 ~~{\rm{and}}~~t>t_0,
\end{align*}
which together with \eqref{eq4-1} deduces that
\begin{align*}
\|v(\cdot,t)\|_{{L^{\infty}}(\Om)}\rightarrow 0 ~~{\rm{as}}~t\rightarrow \infty.
\end{align*}
Consequently, \eqref{con} is proved. Finally, \eqref{th-de3} is a consequence of \eqref{de3}, and \eqref{th-de4} can be derived from \eqref{eq-wabe} and \eqref{eq-abe}. $\hfill\blacksquare$

\textbf{Acknowledgements.}
This work is supported by the Deutsche Forschungsgemeinschaft (No. 462888149).

\small{
}

\end{document}